\documentclass[12pt]{article}

\usepackage{amsmath,amssymb,amsthm}
\usepackage{geometry}
\usepackage{enumitem}
\usepackage{titlesec}
\usepackage{graphicx}
\usepackage{float}
\usepackage{placeins}

\title{Numerical Investigation of Soft-to-Hard Convergence:
Discretization, Admissibility, and Spectral Regularity}

\author{
Sung-Sik Kwon\\
Department of Mathematics and Physics\\
North Carolina Central University\\
\texttt{skwon@nccu.edu}
\and
Marcus Nunda\\
Department of Mathematics and Physics\\
North Carolina Central University\\
\texttt{mnunda@eagles.nccu.edu}
}
\date{}

\begin{document}

\maketitle

\begin{abstract}
This paper presents a numerical investigation of soft-to-hard terminal
constraint convergence for optimal control problems governed by a second-order
ordinary differential equation and the one-dimensional heat equation. The study
complements the analytical convergence theory through comparisons with exact
solutions and independently computed numerical approximations. For the heat
equation, a fully discrete Crank--Nicolson finite-difference formulation is
validated against exact modal solutions, exhibiting approximately second-order
convergence in the terminal-state and midpoint-sampled control errors under
simultaneous space-time mesh refinement. Numerical experiments
demonstrate the predicted first-order penalty convergence for benchmark
problems and investigate the admissibility threshold governing finite-energy
exact terminal controllability. Additional modal experiments examine the
dependence of terminal-error convergence rates on spectral regularity,
illustrating intermediate rates consistent with the theoretical estimates and
asymptotic predictions. The results provide computational evidence supporting
the soft-to-hard convergence framework and clarify the roles of discretization,
target admissibility, and spectral regularity.
\end{abstract}




\section{Introduction}

Terminal state constraints arise naturally in optimal control problems where
the objective is to steer the state of a dynamical system to a prescribed
target at a specified final time. Such problems occur in numerous applications,
including trajectory planning, thermal regulation, inverse problems, and
distributed parameter systems. The mathematical theory of terminally
constrained optimal control has been extensively developed for both
finite-dimensional and infinite-dimensional systems; see, for example,
\cite{Bertsekas1982,Hinze2009,Lions1971,Troltzsch2010}.

Although exact terminal constraints are natural from a modeling perspective,
they often present significant analytical and computational challenges. A
common approach is to replace the hard terminal constraint by a quadratic
penalty on the terminal error, thereby converting the constrained optimization
problem into an unconstrained one. Penalty and regularization methods have long
played an important role in optimization and optimal control because they
simplify both the theoretical analysis and the numerical solution of
constrained problems. Classical references include
\cite{BonnansShapiro2000,Glowinski1984,ItoKunisch2008, RockafellarWets1998},
while applications to PDE-constrained optimization are discussed in
\cite{Casas2012,HintermuellerHinze2006,Maurer2008,MeidnerVexler2008a}.

As the penalty parameter increases, the solution of the penalized problem is
expected to approach the solution of the corresponding hard-constrained
problem, provided the prescribed terminal state is reachable. This phenomenon,
referred to as
\emph{soft-to-hard convergence}, is fundamental to the
mathematical justification of penalty methods and to the development of
numerical algorithms for constrained optimal control.

In our previous work \cite{Kwon2026SoftToHard}, we established a quantitative
soft-to-hard convergence theory for terminally penalized optimal control
problems in both finite-dimensional and infinite-dimensional settings. We
proved that, for admissible terminal targets, the penalized optimal controls
converge strongly to the minimum-energy hard-constrained controls. Under
additional spectral regularity assumptions, we obtained terminal-state
convergence estimates of the form
\[
\|y_\alpha(T)-y_T\|_{L^2(0,1)}
=O(\alpha^{-(1/2+\theta)}),
\qquad 0\le\theta\le\frac12,
\]
and optimal-control convergence estimates of the form
\[
\|u_\alpha^*-u^*\|_{L^2(Q)}
=O(\alpha^{-\theta}),
\qquad 0<\theta\le1.
\]

The corresponding spectral regularity conditions become stronger as the
convergence exponents increase.

The analysis also provides a spectral characterization of admissible terminal
targets for the heat equation, identifying precisely those targets that can be
reached by finite-energy controls.

The purpose of the present paper is to provide a numerical investigation of the
theoretical results developed in
\cite{Kwon2026SoftToHard}. We examine how the predicted
convergence behavior is reflected in numerical computations, investigate the
accuracy of independently computed approximations, and explore the influence of
spectral regularity on terminal accuracy and control energy. Particular
attention is given to the dependence of terminal-error convergence rates on the
regularity of the prescribed target.

The principal contributions of this paper are threefold. First, we provide
numerical verification of the predicted first-order soft-to-hard convergence
behavior for benchmark problems satisfying the corresponding regularity
assumptions. These experiments include comparisons between exact analytical
solutions and independently computed Crank--Nicolson finite-difference
approximations. Second, we investigate the admissibility criterion developed in
\cite{Kwon2026SoftToHard}, illustrating how the spectral decay of terminal data
determines the existence of finite-energy hard-constrained controls and
influences the required control energy. Third, we examine regularity-dependent
terminal-error convergence rates using multimodal targets with power-law
Fourier coefficients. These experiments illustrate intermediate convergence
rates and their transition toward first-order behavior as the spectral
regularity increases.

To achieve these objectives, we consider two complementary model problems that
illustrate the transition from an effectively finite-dimensional prototype to
an infinite-dimensional distributed parameter system.

\begin{enumerate}

\item A second-order ordinary differential equation,
referred to as the \emph{rocket prototype}. Although
the control is sought in the infinite-dimensional
space \(L^2(0,T)\), the terminal condition reduces
to a single scalar moment constraint, giving the
problem an effectively one-dimensional structure.
Because this model admits explicit analytical
solutions, the numerical approximations can be
compared directly with their closed-form counterparts,
providing a convenient benchmark for verifying
soft-to-hard convergence.

\item A distributed optimal control problem governed
by the one-dimensional heat equation. In contrast
to the rocket prototype, the terminal condition
involves infinitely many modal components, leading
to a genuinely infinite-dimensional problem.
This model is studied using modal representations
and fully discrete Crank--Nicolson space--time
approximations. Comparisons between these approaches
provide an independent numerical check of the
predicted convergence behavior for the single-mode
benchmark. Additional modal experiments investigate
admissibility and regularity-dependent convergence
for targets involving multiple Fourier modes.

\end{enumerate}

The remainder of the paper is organized as follows. Section~2 investigates the
rocket prototype, comparing analytical solutions with independently computed
time-discrete approximations. Section~3 presents the heat-equation experiments,
including Crank--Nicolson validation, soft-to-hard convergence, admissibility
thresholds, and intermediate convergence rates associated with spectral
regularity. Section~4 summarizes the numerical findings and discusses their
implications.

\section{Rocket Prototype Problem}
\label{sec:rocket}

The rocket prototype provides a simple computational benchmark for
investigating soft-to-hard convergence before considering the heat equation.
Since the problem admits explicit analytical solutions, it allows direct
comparison with independently computed time-discrete approximations. We use
these comparisons to validate the numerical implementation and examine the
convergence of penalized controls toward the corresponding minimum-energy
hard-constrained controls as the penalty parameter increases.

\subsection{Model Problem}
\label{subsec:rocket_model}

We consider the prototype optimal control problem
\[
y''(t)=v(t)-1,
\qquad 0<t<T,
\]
subject to
\[
y(0)=0,
\qquad
y'(0)=0.
\]
Here, \(y(t)\) denotes the position of the rocket, \(v(t)\) is the control
input, and the constant forcing term represents a normalized gravitational
acceleration.

The hard terminal requirement is
\[
y(T)=y_T.
\]
The corresponding minimum-energy problem is
\[
\min_{v\in L^2(0,T)}
\frac12\int_0^T |v(t)|^2\,dt
\]
subject to the state equation, the initial conditions, and the terminal
constraint \(y(T)=y_T\).

For a penalty parameter \(\alpha>0\), the terminal constraint is replaced by a
quadratic penalty, yielding
\[
\min_{v\in L^2(0,T)}
J_\alpha(v),
\]
where
\[
J_\alpha(v)
=
\frac12\int_0^T |v(t)|^2\,dt
+
\frac{\alpha}{2}
|y(T)-y_T|^2.
\]
The first term measures the control energy, while the second penalizes
deviations from the prescribed terminal target. The theory developed in
\cite{Kwon2026SoftToHard} shows that the penalized solution converges to
the minimum-energy hard-constrained solution as
\(\alpha\to\infty\), with an \(O(\alpha^{-1})\) convergence rate under
the applicable assumptions.


\subsection{Time-Discrete Optimization}
\label{subsec:rocket_validation}
We construct a time-discrete optimization problem by approximating the control
with piecewise-constant functions and integrating the state equations over each
time interval. The resulting finite-dimensional penalized and hard-constrained
problems are solved independently, allowing their convergence behavior to be
examined directly.

To construct the time-discrete problem, we rewrite the second-order state
equation as the first-order system
\[
y'(t)=z(t),
\qquad
z'(t)=v(t)-1,
\]
with initial conditions
\[
y(0)=0,
\qquad
z(0)=0,
\]
where \(z(t)\) denotes the rocket velocity.

Let
\[
t_n=n\Delta t,
\qquad
n=0,\ldots,N,
\qquad
\Delta t=\frac{T}{N},
\]
and approximate the control by a piecewise-constant function,
\[
v(t)=v_n,
\qquad
t_n\le t<t_{n+1},
\qquad
n=0,\ldots,N-1.
\]
Since the control is constant on each interval,
\[
z'(t)=v_n-1,
\qquad
t_n\le t<t_{n+1},
\]
and exact integration over \([t_n,t_{n+1}]\) gives
\[
z_{n+1}
=
z_n+\Delta t\,(v_n-1).
\]
Using
\[
y'(t)=z(t)
=
z_n+(t-t_n)(v_n-1),
\]
integration over the same interval yields
\[
y_{n+1}
=
y_n+\Delta t\,z_n
+\frac{\Delta t^2}{2}(v_n-1),
\qquad
n=0,\ldots,N-1,
\]
with
\[
y_0=0,
\qquad
z_0=0.
\]

Introduce the discrete control vector
\[
\mathbf v=(v_0,\ldots,v_{N-1})^T\in\mathbb R^N.
\]
Applying the recursion repeatedly yields
\[
z_n
=
\Delta t
\sum_{k=0}^{n-1}(v_k-1),
\qquad
n=1,\ldots,N.
\]
Substituting this expression into the recursion for \(y_n\) gives
\[
\begin{aligned}
y_n
&=
\sum_{j=0}^{n-1}
\left[
\Delta t\,z_j
+\frac{\Delta t^2}{2}(v_j-1)
\right] \\
&=
\Delta t^2
\sum_{j=0}^{n-1}
\sum_{k=0}^{j-1}(v_k-1)
+
\frac{\Delta t^2}{2}
\sum_{j=0}^{n-1}(v_j-1),
\qquad
n=1,\ldots,N.
\end{aligned}
\]
In particular, for \(n=N\),
\[
y_N
=
\Delta t^2
\sum_{j=0}^{N-1}
\sum_{k=0}^{j-1}(v_k-1)
+
\frac{\Delta t^2}{2}
\sum_{j=0}^{N-1}(v_j-1).
\]
Observing that each term \(v_k-1\) appears \(N-1-k\) times in the double sum,
and using \(T=N\Delta t\) and \(t_k=k\Delta t\), we obtain
\[
y_N
=
\sum_{k=0}^{N-1}
\Delta t
\left(
T-t_k-\frac{\Delta t}{2}
\right)v_k
-
\frac{T^2}{2}.
\]

Define
\[
w_n
=
\Delta t
\left(
T-t_n-\frac{\Delta t}{2}
\right),
\qquad
\mathbf w=(w_0,\ldots,w_{N-1})^T,
\]
and let
\[
d=y_T+\frac{T^2}{2}.
\]
Then the discrete terminal constraint is simply
\[
\mathbf w^T\mathbf v=d.
\]

Because the control is piecewise constant,
\[
\frac12\int_0^T|v(t)|^2\,dt
=
\frac{\Delta t}{2}\|\mathbf v\|^2.
\]
Hence the penalized problem is
\[
\min_{\mathbf v\in\mathbb R^N}
\left\{
\frac{\Delta t}{2}\|\mathbf v\|^2
+
\frac{\alpha}{2}
|\mathbf w^T\mathbf v-d|^2
\right\}.
\]
The first-order optimality condition is
\[
\left(
\Delta t\,I+\alpha\mathbf w\mathbf w^T
\right)\mathbf v_{\alpha}
=
\alpha d\,\mathbf w,
\]
where \(I\) denotes the \(N\times N\) identity matrix.

The corresponding hard-constrained solution is obtained from
\[
\min_{\mathbf v\in\mathbb R^N}
\frac{\Delta t}{2}\|\mathbf v\|^2
\]
subject to
\[
\mathbf w^T\mathbf v=d.
\]
The associated Karush--Kuhn--Tucker system is
\[
\begin{pmatrix}
\Delta t\,I & \mathbf w\\
\mathbf w^T & 0
\end{pmatrix}
\begin{pmatrix}
\mathbf v^*\\
\lambda
\end{pmatrix}
=
\begin{pmatrix}
\mathbf0\\
d
\end{pmatrix}.
\]

After computing either discrete control, the corresponding position and
velocity are recovered by the state recursions. Thus, both the penalized and
hard-constrained solutions are obtained independently from the discrete
optimization problems rather than by sampling the continuous analytical
solutions.

Because explicit analytical solutions are available for the rocket prototype,
the implementation was verified by comparing the discrete penalized state and
control with the corresponding analytical penalized solution for a
representative value of the penalty parameter. The close agreement confirms the
correctness of the time-discrete optimization procedure. The corresponding
comparison is included in Appendix~\ref{appendix:rocket_validation}. We now
proceed to the soft-to-hard convergence study based entirely on independently
computed discrete penalized and hard-constrained solutions.

\subsection{Discrete Soft-to-Hard Convergence}
\label{subsec:rocket_convergence}

Having verified the time-discrete implementation against the available
closed-form solutions in Section~\ref{subsec:rocket_validation}, we now
investigate the soft-to-hard convergence of the time-discrete optimization
problem. Throughout this experiment, the temporal discretization is held fixed,
with
\[
T=2,
\qquad
y_T=1,
\qquad
N=80,
\qquad
\Delta t=\frac{T}{N}.
\]

The penalized problem is solved for ten logarithmically spaced values
\[
\alpha_k=10^{4(k-1)/9},
\qquad k=1,\ldots,10,
\]
covering the range
\[
10^0\le\alpha\le10^4.
\]
The logarithmic spacing samples four orders of magnitude uniformly on the
logarithmic scale, allowing both pre-asymptotic and asymptotic behavior to be
examined using a modest number of penalty parameters.

For each value of \(\alpha\), the discrete penalized solution
\[
(\mathbf y_{\alpha},\mathbf v_{\alpha})
\]
is compared with the independently computed discrete hard-constrained solution
\[
(\mathbf y^*,\mathbf v^*),
\]
which is obtained by solving the discrete KKT system introduced in
Section~\ref{subsec:rocket_validation}. Thus, the convergence study is
performed entirely within the validated time-discrete framework.

Three errors are monitored. The terminal error is
\[
E_T(\alpha)
=
|y_{\alpha}(T)-y_T|.
\]
Since the discrete hard-constrained solution satisfies
\[
y^*(T)=y_T
\]
up to numerical roundoff, this is also the difference between the terminal
values of the penalized and hard-constrained discrete states.

The discrete control error is
\[
E_v(\alpha)
=
\|\mathbf v_{\alpha}-\mathbf v^*\|,
\]
where the discrete norm is defined by
\[
\|\mathbf q\|
=
\left(
\Delta t
\sum_{n=0}^{N-1}|q_n|^2
\right)^{1/2}.
\]
The control vector is naturally interpreted as a piecewise-constant function on
the time intervals, so this norm is the corresponding discrete \(L^2(0,T)\)
norm.

The discrete state error is
\[
E_y(\alpha)
=
\|\mathbf y_{\alpha}-\mathbf y^*\|_{\mathrm{trap}},
\]
where the state vectors contain the nodal values
\[
\mathbf y=(y_0,\ldots,y_N)^\top,
\]
and the trapezoidal discrete norm is defined by
\[
\|\mathbf q\|_{\mathrm{trap}}
=
\left[
\Delta t
\left(
\frac12|q_0|^2
+
\sum_{n=1}^{N-1}|q_n|^2
+
\frac12|q_N|^2
\right)
\right]^{1/2}.
\]
This norm approximates the continuous \(L^2(0,T)\) norm of the state difference
using the computed nodal values.

For an error of the form
\[
E(\alpha)\approx C\alpha^p,
\]
the exponent \(p\) is estimated by a least-squares fit of
\(\log_{10}E(\alpha)\) against \(\log_{10}\alpha\). Slopes are computed
both over the full range of penalty parameters and over the asymptotic range
\[
\alpha\geq10^2.
\]
Restricting the fit to the larger penalty parameters reduces the influence of
the pre-asymptotic regime and provides a more reliable estimate of the
asymptotic convergence rate.

\begin{figure}[H]
\centering
\includegraphics[width=0.85\textwidth]
{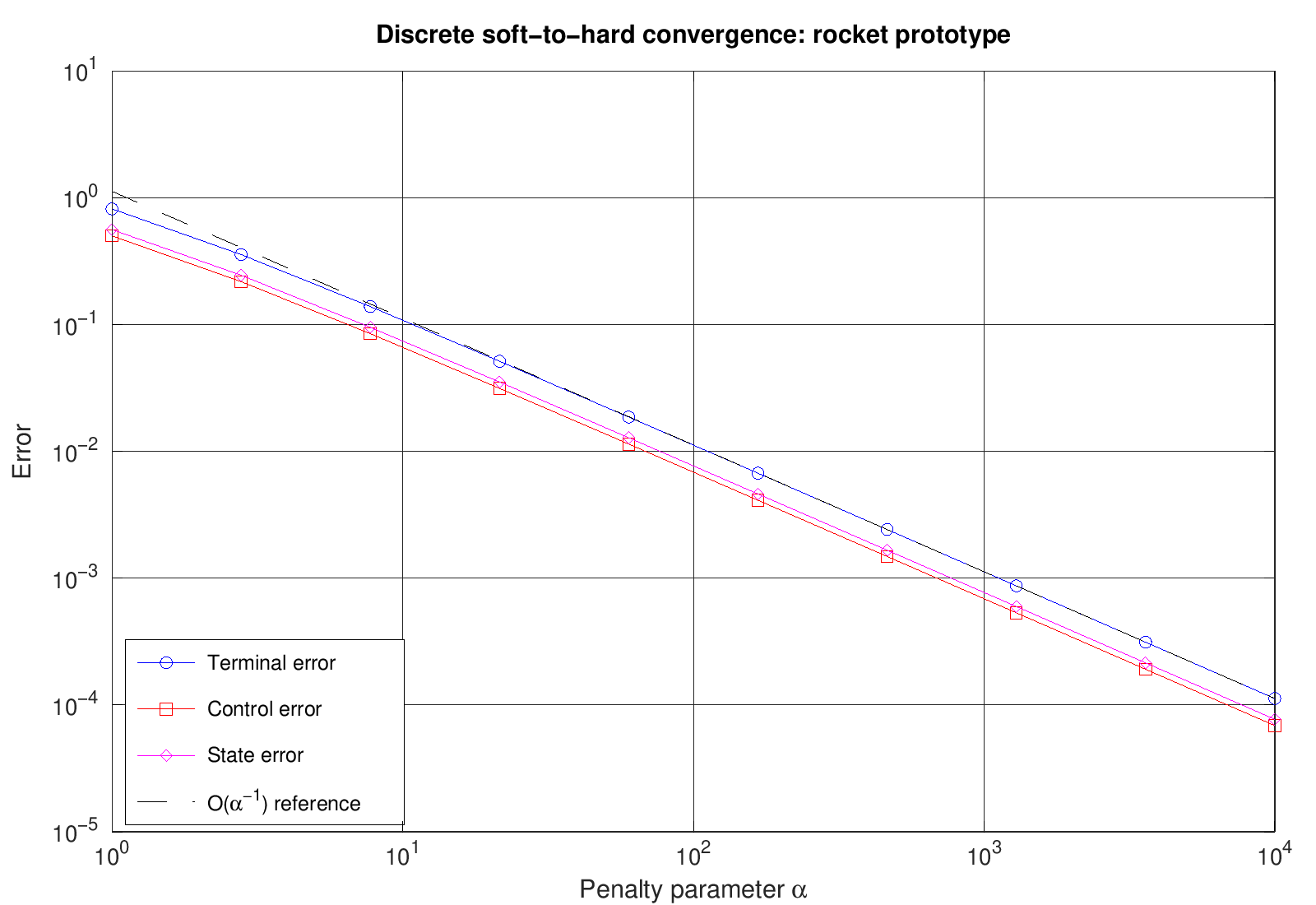}
\caption{Discrete soft-to-hard convergence for the rocket
prototype. Terminal, control, and state errors are plotted
against the penalty parameter \(\alpha\) on logarithmic
scales. The reference line has slope \(-1\).}
\label{fig:rocket_convergence}
\end{figure}

Figure~\ref{fig:rocket_convergence} shows the terminal, control, and state
errors together with a reference line of slope \(-1\).

The least-squares slopes obtained from the three error curves are reported in
Table~\ref{tab:rocket_slopes}.

\begin{table}[H]
\centering
\caption{Least-squares slopes of the discrete error curves.}
\label{tab:rocket_slopes}
\begin{tabular}{lcc}
\hline
Error & Full range & \(\alpha\ge10^2\) \\
\hline
Terminal error \(E_T\) & \(-0.976083\) & \(-0.999500\) \\
Control error \(E_v\)  & \(-0.976083\) & \(-0.999500\) \\
State error \(E_y\)    & \(-0.976083\) & \(-0.999500\) \\
\hline
\end{tabular}
\end{table}

As shown in Figure~\ref{fig:rocket_convergence}, the three error curves
approach the reference slope
\(-1\) as \(\alpha\) increases. Over the full parameter
range, the fitted slope is \(-0.976083\) for each error. Restricting the
least-squares fit to \(\alpha\ge10^2\) gives the asymptotic slope
\(-0.999500\), which is very close to the theoretical value \(-1\). These
results are consistent with
\[
E_T(\alpha)=O(\alpha^{-1}),\qquad
E_v(\alpha)=O(\alpha^{-1}),\qquad
E_y(\alpha)=O(\alpha^{-1}).
\]

The equality of the three fitted slopes is a special feature of the scalar
rocket prototype. Indeed, the discrete optimality systems imply
\[
\mathbf v_\alpha-\mathbf v^*
=
-\frac{1}{1+\alpha a_{\Delta t}}\mathbf v^*,
\qquad
a_{\Delta t}=\frac{\mathbf w^\top\mathbf w}{\Delta t}.
\]
Since the discrete state depends affinely on the control, the terminal,
control, and state errors are all proportional to the common factor
\[
\frac{1}{1+\alpha a_{\Delta t}}.
\]
Consequently, their log--log curves differ only by constant vertical shifts and
have identical slopes. As \(\alpha\to\infty\), this factor behaves like
\(O(\alpha^{-1})\), explaining the observed first-order convergence.

Thus, for a fixed temporal discretization, the discrete penalized solutions
reproduce the soft-to-hard convergence rate predicted by the continuous theory.
The deviation of the full-range slope from \(-1\) is attributable to smaller
penalty parameters, for which the asymptotic behavior is not yet fully
apparent.

This experiment should be interpreted as a fixed-mesh soft-to-hard convergence
study. It establishes convergence of the discrete penalized solution toward the
discrete hard-constrained solution as the penalty parameter increases, but does
not measure convergence with respect to the time-step size.

The rocket prototype therefore confirms that the discrete penalized solutions
exhibit the
\(O(\alpha^{-1})\) soft-to-hard convergence
behavior predicted by the continuous theory in \cite{Kwon2026SoftToHard}.
Having established this convergence behavior in the time-discrete setting, we
next turn to the heat equation.


\section{Heat Equation Numerical Experiments}

We now investigate the numerical behavior of terminal penalization for the
one-dimensional heat equation, focusing on the accuracy of the numerical
approximations and their consistency with the theoretical soft-to-hard
convergence results established in
\cite{Kwon2026SoftToHard}.


\subsection{Numerical Framework}

The numerical experiments employ two complementary computational approaches.
The first is a finite-mode truncation of the explicit modal representation
developed in \cite{Kwon2026SoftToHard}, which serves as a reference solution.
The second is an independently implemented finite-difference approximation
based on centered spatial differences and Crank--Nicolson time stepping.

Throughout the heat-equation numerical experiments, we take
\(y_0=0\) and impose homogeneous Dirichlet boundary conditions.
This choice simplifies the numerical formulation while retaining the principal
convergence and admissibility phenomena investigated in this paper.

In the modal computation, the terminal state is approximated by the truncated
eigenfunction expansion
\[
y_\alpha(x,T)
=
\sum_{n=1}^{N}y_{\alpha,n}(T)e_n(x),
\]
where
\[
e_n(x)=\sqrt{2}\sin(n\pi x),
\qquad
\lambda_n=(n\pi)^2.
\]

For each retained mode, the penalized terminal coefficient is given by
\[
y_{\alpha,n}(T)
=
e^{-\lambda_nT}y_{0,n}
+
\frac{\alpha a_n}{1+\alpha a_n}d_n,
\]
where
\[
d_n=y_{T,n}-e^{-\lambda_nT}y_{0,n},
\qquad
a_n=\frac{1-e^{-2\lambda_nT}}{2\lambda_n}.
\]
Under the zero-initial-state assumption, this reduces to
\[
y_{\alpha,n}(T)
=
\frac{\alpha a_n}{1+\alpha a_n}d_n,
\qquad d_n=y_{T,n}.
\]

For the fully discrete approximation, the spatial interval
\((0,L)\) is divided into \(m\) subintervals of length
\(h=L/m\), and the time interval \((0,T)\) is divided
into \(n\) subintervals of length \(k=T/n\). Let
\[
\mathbf y^j
=
\begin{bmatrix}
U_1^j & \cdots & U_{m-1}^j
\end{bmatrix}^{\top}
\]
denote the vector of interior state values at time \(t_j\), and let
\(\mathbf v^j\) denote the intervalwise-constant control on the
\(j\)-th time interval. The Crank--Nicolson equations take the form

\[
B\mathbf y^j-C\mathbf y^{j-1}
=
k\mathbf v^j,
\qquad
j=1,\ldots,n,
\]
where
\[
B
=
I-\frac{k}{2}L_h,
\qquad
C
=
I+\frac{k}{2}L_h,
\]
and \(L_h\) is the centered finite-difference approximation of the Dirichlet
Laplacian.

Stacking the state and control vectors as
\[
\mathbf y
=
\begin{bmatrix}
\mathbf y^1\\
\vdots\\
\mathbf y^n
\end{bmatrix},
\qquad
\mathbf v
=
\begin{bmatrix}
\mathbf v^1\\
\vdots\\
\mathbf v^n
\end{bmatrix},
\]
gives the all-at-once space--time system
\[
A\mathbf y
=
k\mathbf v+\mathbf b_0.
\]

Here \(\mathbf b_0\) contains the contribution of the prescribed initial
condition. In the experiments considered below, the initial state is zero, so
\(\mathbf b_0=\mathbf0\). The terminal state is included directly in \(\mathbf
y\), and the extraction matrix \(M\) satisfies
\[
M\mathbf y=\mathbf y^n.
\]

The discrete penalized optimal control problem is
\[
\min_{\mathbf y,\mathbf v}
\left\{
\frac{hk}{2}\|\mathbf v\|^{2}
+
\frac{\alpha h}{2}
\|M\mathbf y-\mathbf y_T\|^{2}
\right\},
\]
subject to
\[
A\mathbf y=k\mathbf v,
\]
where \(\|\cdot\|\) denotes the Euclidean norm.

Eliminating the control via
\[
\mathbf v=\frac{1}{k}A\mathbf y
\]
yields the reduced optimization problem
\[
\min_{\mathbf y}
\left\{
\frac{h}{2k}\|A\mathbf y\|^{2}
+
\frac{\alpha h}{2}
\|M\mathbf y-\mathbf y_T\|^{2}
\right\}.
\]

For each penalty parameter \(\alpha\), the optimal discrete state is obtained
by solving the symmetric positive-definite linear system
\[
\left(
\frac{h}{k}A^{\top}A
+
\alpha hM^{\top}M
\right)
\mathbf y_\alpha
=
\alpha hM^{\top}\mathbf y_T.
\]

The corresponding discrete control and terminal state are then recovered from
\[
\mathbf v_\alpha
=
\frac{1}{k}A\mathbf y_\alpha,
\qquad
\mathbf y_\alpha^{\,n}
=
M\mathbf y_\alpha.
\]

The terminal error and control norm reported below are defined by
\[
E_T(\alpha)
=
\left(
h
\sum_{i=1}^{m-1}
\left|
(\mathbf y_\alpha^{\,n})_i-(\mathbf y_T)_i
\right|^2
\right)^{1/2},
\]
and
\[
\|\mathbf v_\alpha\|
=
\left(
hk
\sum_{j=1}^{n}
\sum_{i=1}^{m-1}
|v_{\alpha,i}^{\,j}|^2
\right)^{1/2},
\]
which are finite-difference approximations of the continuous
\(L^2(0,L)\) terminal error and \(L^2((0,L)\times(0,T))\) control norm,
respectively.

The modal computation provides a reference solution, while the fully discrete
formulation independently verifies that the predicted soft-to-hard convergence
is reproduced by the space--time discretization.


\subsection{Validation of the Fully Discrete Approximation}

To validate the fully discrete Crank--Nicolson implementation, we compare the
finite-difference solution with the exact penalized modal solution for the
single-mode terminal target
\[
y_T(x)=\sin(\pi x),
\]
with zero initial state. Since only the first eigenmode is active, the exact
penalized solution is available in closed form. This comparison isolates the
discretization error at a fixed penalty parameter, independently of the
soft-to-hard convergence investigated in the following subsection.

The penalty parameter is fixed at
\[
\alpha=10^4,
\]
while the spatial and temporal meshes are simultaneously refined according to
\[
m=10,\;20,\;40,\;80,
\qquad
n=2m.
\]

For each mesh, the computed terminal state and control are compared with the
corresponding exact modal solutions evaluated on the spatial and space-time
grids. The terminal-state and control errors are defined by
\[
E_y
=
\sqrt{h}\,
\left\|
\mathbf y^{\,n}-\mathbf y(T)
\right\|,
\]
and
\[
E_v
=
\sqrt{hk}\,
\left\|
\mathbf v-\mathbf v^{\mathrm{exact}}
\right\|,
\]
where \(\mathbf y^{\,n}\) and \(\mathbf v\) denote the computed
finite-difference solutions, while \(\mathbf y(T)\) and
\(\mathbf v^{\mathrm{exact}}\) denote the exact penalized modal
solutions evaluated on the corresponding grids. The vector norms in these
expressions are Euclidean norms, so the errors represent weighted discrete
\(L^2\) norms.

Since the fully discrete control is represented by one value on each time
interval, the exact penalized control is evaluated at the temporal midpoints
\[
t_{j-\frac12}=\left(j-\frac12\right)k,
\qquad j=1,\ldots,n,
\]
when computing the control error.

The mesh-refinement results, reported in Table~\ref{tab:mesh-refinement} in
Appendix~\ref{appendix:heat_validation}, show that both the terminal-state
error and the control error exhibit approximately second-order convergence
under simultaneous space-time refinement, with observed convergence orders
approaching
\(2\) as the mesh is refined. These results provide numerical
validation of the fully discrete Crank--Nicolson implementation for the
single-mode benchmark before it is used to investigate soft-to-hard convergence
in the following subsection.


\subsection{Soft-to-Hard Convergence}

Having validated the fully discrete implementation, we next investigate whether
the numerical solutions reproduce the soft-to-hard convergence predicted by the
continuous theory.

The single-mode terminal target
\[
y_T(x)=\sin(\pi x)
\]
is retained throughout this experiment, with zero initial state. The penalized
problem is solved for representative values of the penalty parameter spanning
four orders of magnitude,
\[
1\le\alpha\le10^4.
\]

\begin{figure}[ht]
\centering
\includegraphics[width=0.75\textwidth]{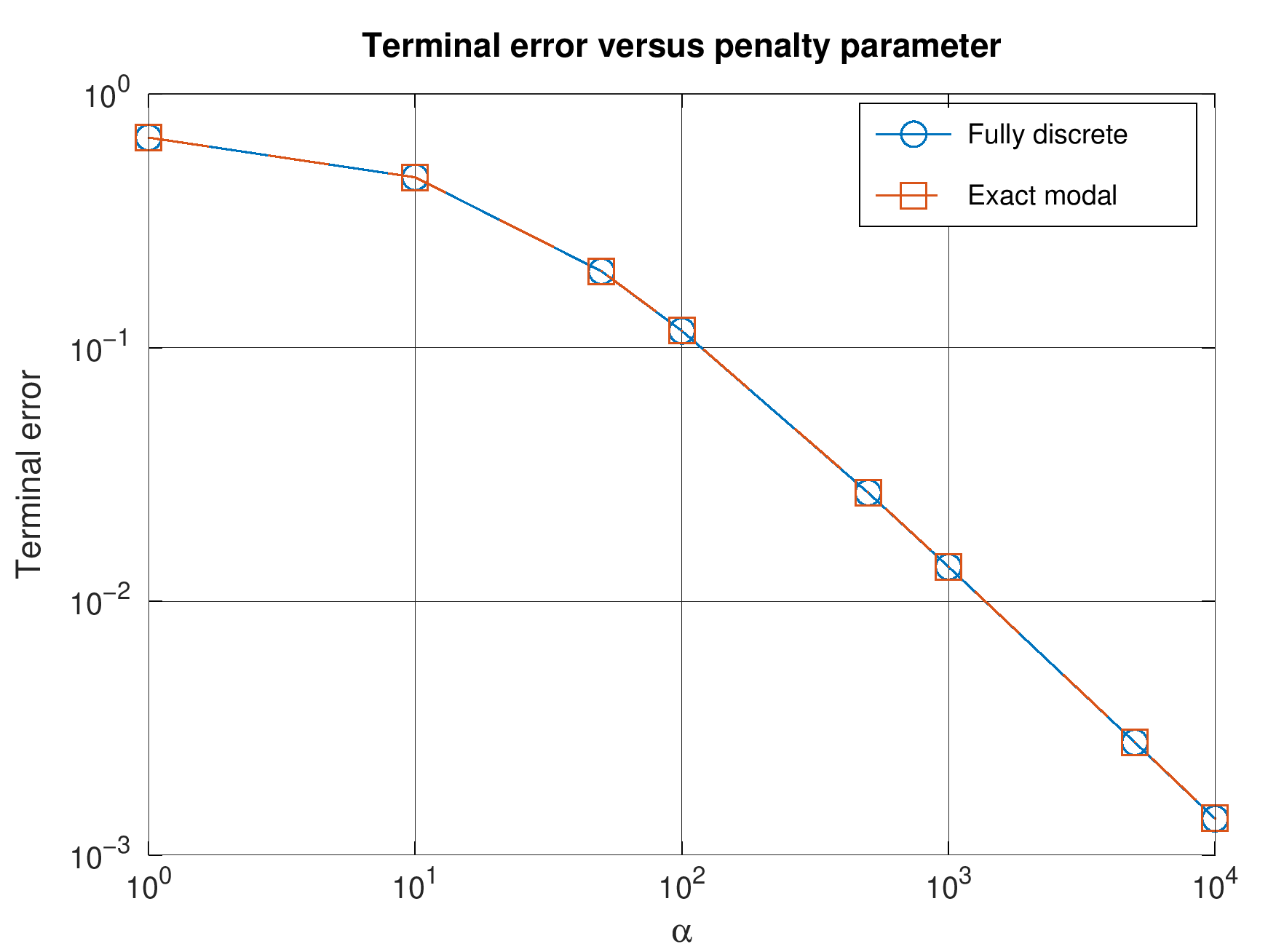}
\caption{Terminal error for the fully discrete and exact modal
solutions as a function of the penalty parameter.}
\label{fig:heat_terminal_error}
\end{figure}

Figure~\ref{fig:heat_terminal_error} shows the terminal error as a function of
the penalty parameter for both the exact modal solution and the fully discrete
Crank--Nicolson approximation. The two curves agree closely over the tested
range, indicating that the fully discrete approximation reproduces the
convergence behavior of the exact modal solution.

For the single-mode target considered here, the spectral regularity assumptions
of the continuous theory are satisfied. Consequently, the continuous penalized
solution obeys the terminal-error estimate
\[
\|y_\alpha(T)-y_T\|_{L^2(0,1)}
=O(\alpha^{-1})
\qquad \text{as } \alpha\to\infty.
\]

As illustrated in Figure~\ref{fig:heat_terminal_error}, both the fully discrete
and exact modal error curves approach a slope of \(-1\) on the log--log scale
as the penalty parameter increases. Using the range
\[
\alpha\ge500,
\]
least-squares fits yield the observed slopes
\[
-0.988333
\quad\text{(fully discrete)},
\qquad
-0.988309
\quad\text{(exact modal)}.
\]
These values agree closely with the theoretical asymptotic slope of \(-1\).

When all penalty parameters are included in the least-squares fits, the
observed slopes are
\[
-0.716814
\quad\text{(fully discrete)},
\qquad
-0.716596
\quad\text{(exact modal)}.
\]
The smaller magnitudes of these slopes reflect the pre-asymptotic behavior
visible in Figure~\ref{fig:heat_terminal_error}. For moderate values of
\(\alpha\), the denominator
\(1+\alpha a_1\) in the exact modal error formula
is not yet dominated by \(\alpha a_1\). Thus, the first-order asymptotic
behavior becomes apparent when \(\alpha a_1\gg1\).

The close agreement between the fully discrete and exact modal results
indicates that the Crank--Nicolson approximation reproduces the predicted
first-order penalty convergence behavior over the tested parameter range for
this single-mode benchmark.


\subsection{Admissibility Threshold Experiments}

The continuous theory developed in \cite{Kwon2026SoftToHard} shows that the
existence of a finite-energy hard-constrained control is determined by the
modal correction coefficients
\[
d_n
=
y_{T,n}
-
e^{-\lambda_nT}y_{0,n},
\]
which combine the prescribed terminal state and the initial condition.

The exact spectral admissibility criterion is
\[
\sum_{n=1}^{\infty}
\frac{|d_n|^2}{a_n}<\infty,
\qquad
a_n=\frac{1-e^{-2\lambda_nT}}{2\lambda_n},
\]
where \(\lambda_n=(n\pi)^2\). Since \(T>0\) is fixed and
\[
\frac{1}{a_n}
=\frac{2\lambda_n}{1-e^{-2\lambda_nT}}
\asymp\lambda_n,
\]
this criterion is equivalent to
\[
\sum_{n=1}^{\infty}
\lambda_n|d_n|^2<\infty.
\]

Since both the optimal hard-constrained control and its penalized counterpart
depend on the initial and terminal states only through the coefficients
\(\{d_n\}\), we prescribe these modal corrections directly in the numerical
experiments. Equivalently, we may take \(y_0=0\), so that
\(d_n=y_{T,n}\) and the prescribed coefficients
represent the Fourier coefficients of the terminal target.

To investigate the influence of regularity, we consider the family of modal
corrections
\[
d_n\propto n^{-p},
\]
where the exponent \(p\) controls the decay of the higher Fourier modes. For
each truncation level \(N\), the coefficients are normalized by
\[
d_n^{(N)}
=
\frac{n^{-p}}
{\left(\sum_{j=1}^{N}j^{-2p}\right)^{1/2}},
\qquad 1\le n\le N.
\]
This normalization gives unit \(\ell^2\)-norm for each truncated coefficient
vector, allowing comparisons across different values of \(p\) without
differences in the overall magnitude of the prescribed modal correction.

Since
\[
\lambda_n|d_n|^2
\asymp
n^{2-2p},
\]
the admissibility condition reduces to the convergence of the series
\[
\sum_{n=1}^{\infty}
n^{2-2p},
\]
which holds if and only if
\[
p>\frac32.
\]

Thus, the critical exponent \(p=\frac32\) separates admissible and
non-admissible terminal data. We consider the representative values
\[
p=1,\quad
1.25,\quad
1.5,\quad
1.75,\quad
2,
\]
which include both non-admissible (\(p\le\frac32\)) and admissible
(\(p>\frac32\)) cases.

Since the hard-constrained control is represented by an infinite modal
expansion, we approximate its norm using the first \(N\) Fourier modes. For
admissible targets, these truncated norms should remain uniformly bounded as
\(N\) increases, whereas for non-admissible targets they are expected to
continue growing.

\begin{figure}[ht]
\centering
\includegraphics[width=0.75\textwidth]{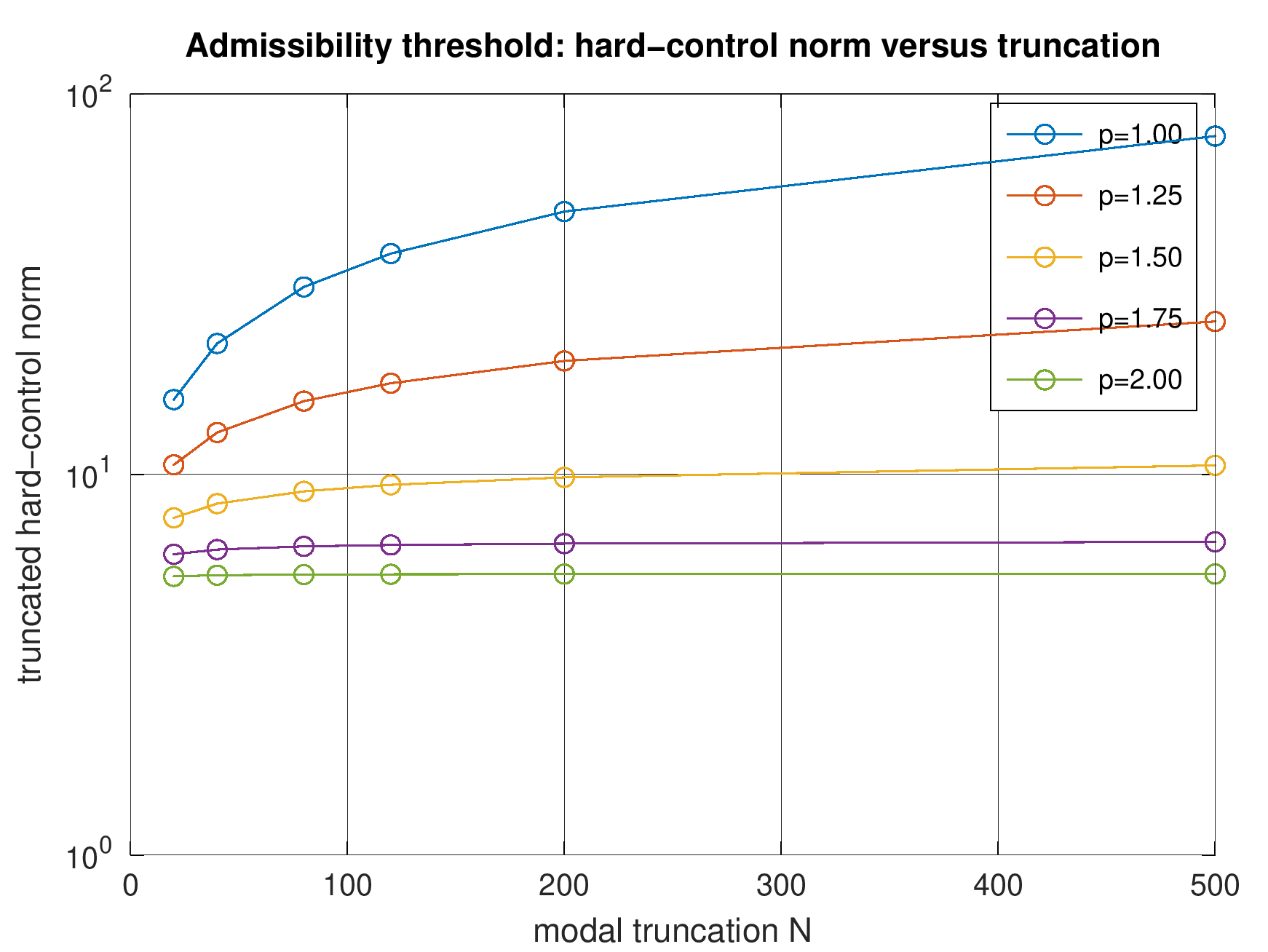}
\caption{Truncated hard-control norms for modal correction
coefficients satisfying \(d_n\propto n^{-p}\). The norms remain bounded
for admissible targets (\(p>\frac32\)) but increase with the number of
retained modes for non-admissible targets (\(p\le\frac32\)).}
\label{fig:hard_control_norms}
\end{figure}

Figure~\ref{fig:hard_control_norms} illustrates this distinction over the
tested truncation levels. For \(p>\frac32\), the truncated hard-control norms
approach finite limiting values, consistent with the existence of a
finite-energy hard-constrained control. In contrast, for \(p\le\frac32\), the
truncated norms continue to increase as additional modes are retained,
consistent with the nonexistence of a finite-energy hard-constrained control in
\(L^2(Q)\) for the corresponding infinite-dimensional problem. Although only
finitely many modes are included in the computation, the observed behavior
agrees with the admissibility criterion established in the continuous theory.

\begin{figure}[ht]
\centering
\includegraphics[width=0.75\textwidth]{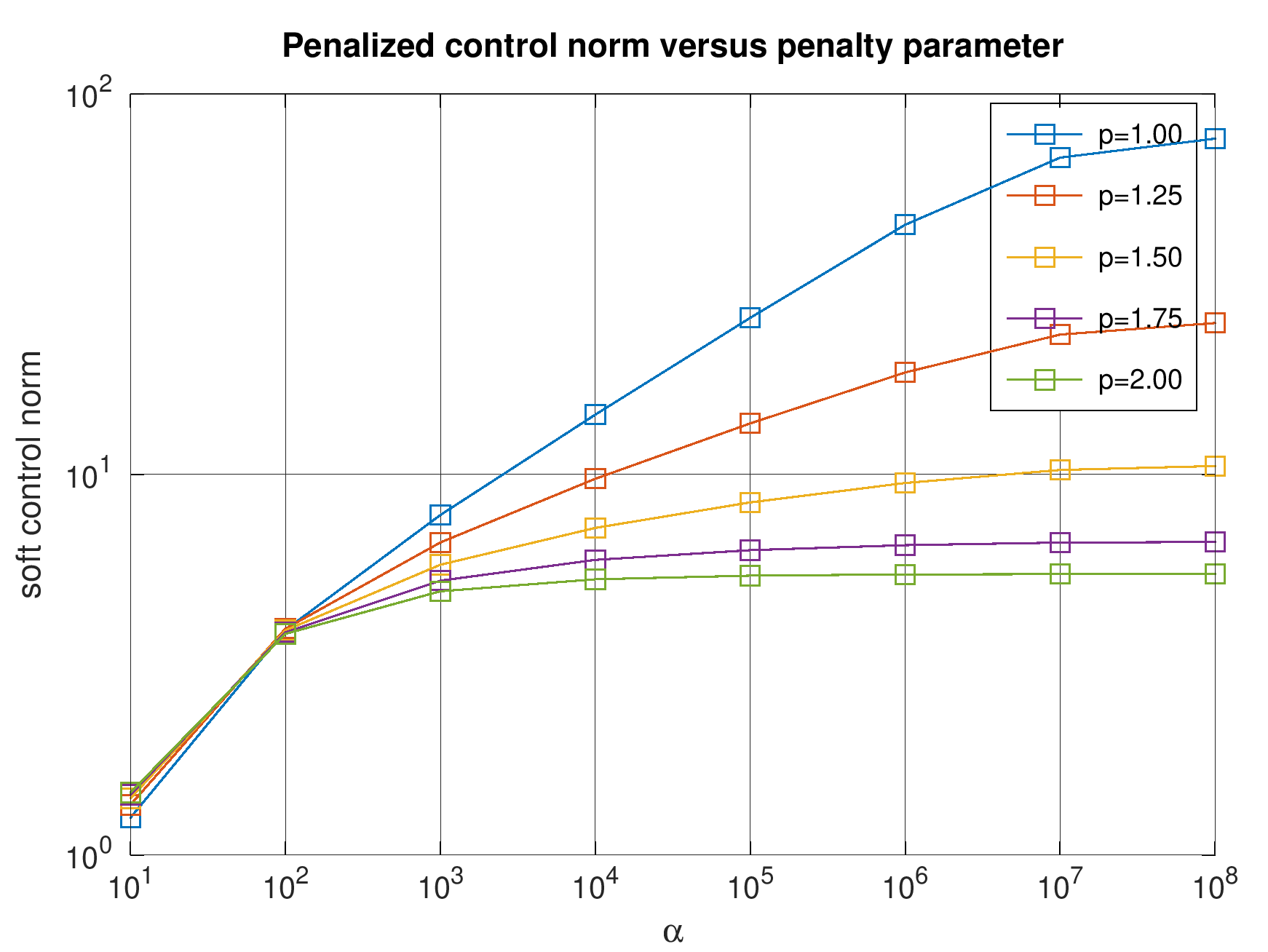}

\caption{Penalized control norms as functions of the penalty
parameter for normalized modal corrections
\(d_n^{(N)}\propto n^{-p}\), with \(N=500\).
The curves illustrate the dependence of control energy
on spectral decay over the tested penalty range.}
\label{fig:admissibility_control_norms}
\end{figure}

The penalized controls exhibit a related qualitative distinction.
Figure~\ref{fig:admissibility_control_norms} displays the control norms as
functions of the penalty parameter for a fixed modal truncation \(N=500\). Over
the computed range of penalty parameters, the control norms associated with
non-admissible targets continue to increase, whereas those associated with
admissible targets show a tendency toward stabilization.

For any fixed modal truncation, however, the penalized control norms approach
finite limits as \(\alpha\to\infty\), regardless of the value of \(p\). The
distinction between admissible and non-admissible targets arises in the
infinite-dimensional limit as the number of retained modes increases. The
observed behavior is consistent with the theoretical admissibility criterion.

These experiments provide numerical illustrations of the admissibility
threshold established by the continuous theory. Although the penalized problem
is well posed for every prescribed modal correction corresponding to an \(L^2\)
terminal target, the computations reveal a qualitative distinction between
admissible and non-admissible data, consistent with the analysis of
\cite{Kwon2026SoftToHard}. Failure of the admissibility condition does not
prevent the penalized terminal states from approximating the prescribed target.
Rather, the admissibility condition governs the existence of a finite-energy
hard-constrained control. In the infinite-dimensional problem, non-admissible
targets require unbounded control energy as the terminal error tends to zero,
whereas admissible targets admit finite-energy hard-constrained controls.


\subsection{Spectral Regularity and Intermediate Convergence Rates}

The preceding experiments examined the admissibility threshold for
finite-energy exact terminal controllability. We now investigate how the
spectral decay of the modal correction coefficients influences the rate of
soft-to-hard convergence.

The continuous theory developed in \cite{Kwon2026SoftToHard} establishes the
terminal-error estimate
\[
\|y_\alpha(T)-y_T\|_{L^2(0,1)}
=O\bigl(\alpha^{-(1/2+\theta)}\bigr),
\qquad 0\le\theta\le\frac12,
\]
provided that the modal corrections satisfy the additional spectral regularity
condition
\[
\sum_{n=1}^{\infty}
\frac{|d_n|^2}{a_n^{1+2\theta}}<\infty.
\]

To examine intermediate convergence rates, we consider normalized modal
corrections satisfying
\[
d_n\propto n^{-p},
\]
with
\[
p=1.75,\quad 2.00,\quad 2.25,\quad 3.00.
\]
All these values satisfy the admissibility condition
\(p>3/2\), but they correspond to different degrees
of spectral regularity.

For the family of modal corrections \(d_n\propto n^{-p}\), the spectral
regularity condition becomes
\[
\sum_{n=1}^{\infty}
n^{2+4\theta-2p}<\infty,
\]
which holds precisely when
\[
p>\frac32+2\theta.
\]
Thus, the admissible range of \(\theta\) depends on the spectral decay exponent
\(p\).

Since \(a_n\asymp n^{-2}\), the dominant contributions to the terminal-error
series can be estimated by splitting the sum near \(n\asymp\sqrt{\alpha}\). For
\(\frac12<p<\frac52\), this gives
\[
\begin{aligned}
E_T(\alpha)^2
&\asymp
\alpha^{-2}
\sum_{n\le\sqrt{\alpha}}n^{4-2p}
+
\sum_{n>\sqrt{\alpha}}n^{-2p}\\
&\asymp
\alpha^{-(p-1/2)}.
\end{aligned}
\]
Consequently,
\[
E_T(\alpha)
\asymp
\alpha^{-(p-1/2)/2},
\qquad \frac12<p<\frac52.
\]

At the borderline case \(p=\frac52\), a logarithmic correction appears:
\[
E_T(\alpha)
\asymp
\alpha^{-1}\sqrt{\log\alpha}.
\]
For \(p>\frac52\), the terminal error instead satisfies
\[
E_T(\alpha)\asymp\alpha^{-1}.
\]

The numerical experiment evaluates the truncated modal error using
\(N=20{,}000\) Fourier modes and 25 logarithmically spaced penalty parameters
over the range
\[
10^2\le\alpha\le10^8.
\]
The observed convergence rates are estimated by least-squares fits of
\(\log_{10}E_T(\alpha)\) against \(\log_{10}\alpha\), using the five largest
penalty parameters.

\begin{table}[ht]
\centering

\caption{Comparison of the predicted infinite-dimensional
asymptotic slopes and observed terminal-error slopes
for modal corrections \(d_n\propto n^{-p}\).
The observed slopes are obtained by least-squares
fits using the five largest penalty parameters
and a modal truncation of \(N=20{,}000\).}
\label{tab:intermediate_rates}

\begin{tabular}{ccc}
\hline
\(p\) & Predicted slope & Observed slope \\
\hline
1.75 & \(-0.625\) & \(-0.62559\) \\
2.00 & \(-0.750\) & \(-0.74999\) \\
2.25 & \(-0.875\) & \(-0.87182\) \\
3.00 & \(-1.000\) & \(-0.99970\) \\
\hline
\end{tabular}
\end{table}

The observed slopes agree closely with the predicted infinite-dimensional
asymptotic rates. In particular, the results demonstrate intermediate
convergence behavior for \(p=1.75\), \(p=2.00\), and \(p=2.25\), while the
faster-decaying coefficients corresponding to \(p=3.00\) yield approximately
first-order convergence.

To assess the influence of modal truncation, the computations were repeated
using
\(N=500,2000,5000,10000,\) and \(20000\).
The fitted slopes approach stable values as
\(N\) increases, with the largest truncation
effects occurring for targets with slower modal decay. At \(N=20000\), the
observed slopes are close to their predicted asymptotic values, supporting the
adequacy of the chosen truncation level over the tested penalty range.

These experiments provide numerical evidence that the terminal-error
convergence rate depends on the spectral regularity of the prescribed target.
For the family of modal corrections considered here, slower spectral decay
leads to intermediate convergence rates, whereas sufficiently rapid decay
yields first-order convergence. The observed rates are consistent with the
regularity-dependent estimates established in \cite{Kwon2026SoftToHard} and
with the more precise asymptotic predictions derived for this family of
targets.

\FloatBarrier

\section{Discussion and Summary}

The numerical experiments provide evidence consistent with the continuous
soft-to-hard convergence theory developed in
\cite{Kwon2026SoftToHard}. The rocket prototype exhibits
first-order convergence of the terminal-state, control, and state-trajectory
errors. For the heat equation, the close agreement between exact modal
solutions and independently computed Crank--Nicolson approximations
demonstrates that the fully discrete method captures the predicted convergence
behavior for the single-mode benchmark.

The experiments also distinguish pre-asymptotic and asymptotic convergence
behavior. As the penalty parameter increases, the observed terminal-error
slopes approach the predicted first-order rate for the benchmark problems
considered.

The admissibility experiments illustrate the critical threshold \(p=3/2\) for
modal corrections
\(d_n\propto n^{-p}\). The truncated hard-control norms
exhibit stabilization for admissible targets and continued growth over the
tested truncation levels for non-admissible targets. The penalized-control
experiments further illustrate how spectral decay influences the control energy
required to approximate the terminal target.

The spectral regularity experiments extend the numerical investigation beyond
first-order convergence. For admissible targets with different modal decay
rates, the observed intermediate terminal-error slopes agree closely with the
asymptotic predictions derived from the modal representation. These results
illustrate the dependence of convergence rates on additional spectral
regularity, consistent with the theoretical estimates.

Overall, the numerical study supports the soft-to-hard convergence framework
through complementary analytical benchmarks, fully discrete validation, and
modal experiments. The results demonstrate the effectiveness of quadratic
terminal penalization for the model problems considered, while distinguishing
the roles of target admissibility, spectral regularity, and numerical
discretization.

\appendix

\section{Validation of the Rocket Prototype}
\label{appendix:rocket_validation}

The rocket prototype admits explicit analytical expressions for both the
penalized state and the corresponding optimal control, making it a natural
benchmark for validating the numerical implementation.

For a representative penalty parameter, the numerical solution is compared with
the analytical solution
\[
y_\alpha(t)
=
\frac12(c_\alpha T-1)t^2
-
\frac{c_\alpha}{6}t^3,
\]
and
\[
v_\alpha(t)
=
c_\alpha(T-t),
\qquad
c_\alpha
=
\frac{\alpha d}{1+\alpha a}.
\]

\begin{figure}[ht]
\centering
\includegraphics[width=0.78\textwidth]{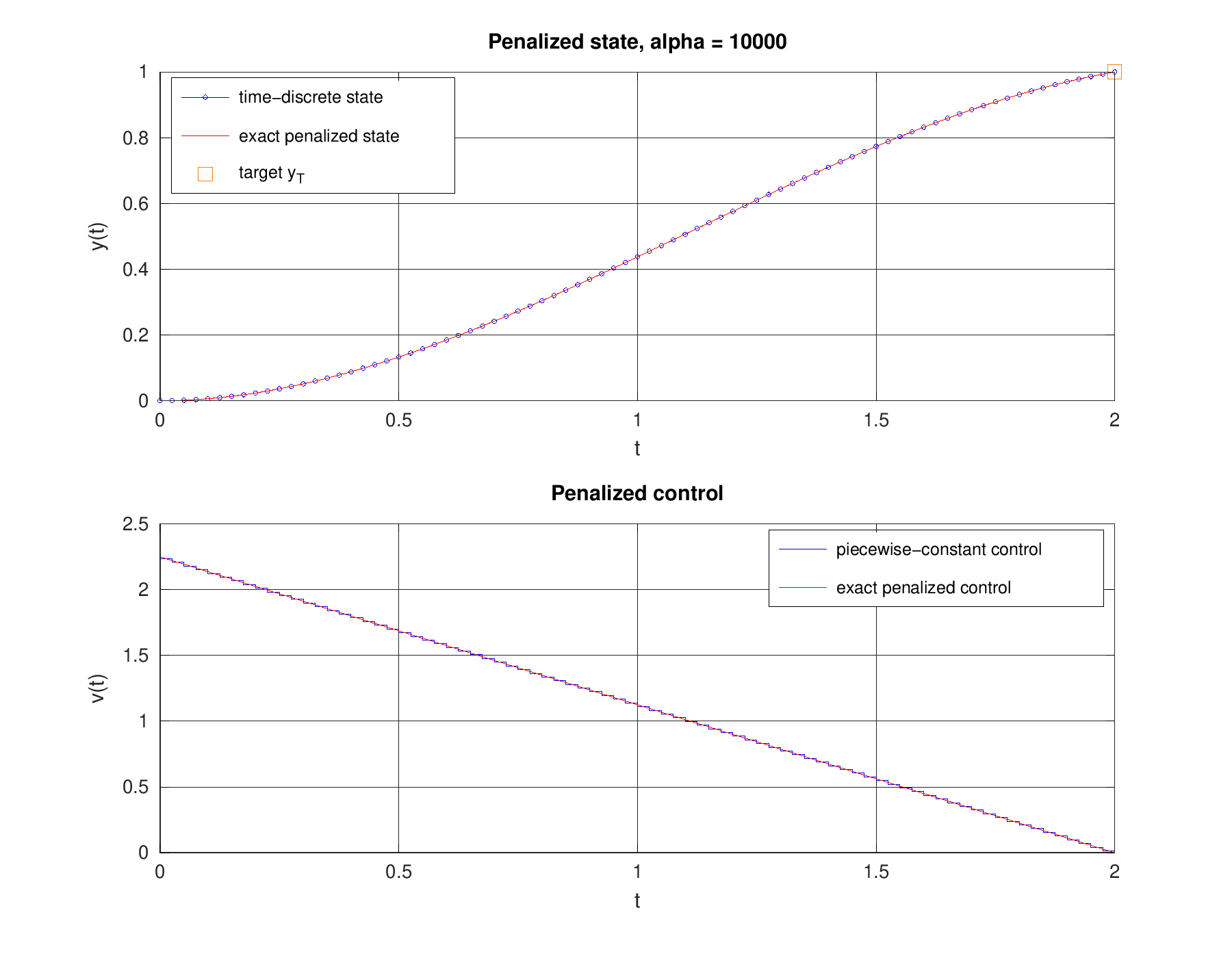}

\caption{Comparison of the time-discrete and exact analytical
penalized solutions for the rocket prototype at
\(\alpha=10^4\). The upper panel compares the position
trajectories, while the lower panel compares the
piecewise-constant discrete control with the continuous
analytical control.}

\label{fig:rocket_validation}
\end{figure}

Figure~\ref{fig:rocket_validation} compares the time-discrete state with the
analytical solution at the discrete time nodes and the piecewise-constant
control with the continuous analytical control. The close agreement provides an
independent numerical validation of the time-discrete implementation.

\FloatBarrier
\section{Validation of the Fully Discrete Heat Equation}
\label{appendix:heat_validation}

The accuracy of the fully discrete Crank--Nicolson implementation is assessed
by comparison with the exact penalized modal solution for the single-mode
terminal target
\[
y_T(x)=\sin(\pi x),
\]
with zero initial state and fixed penalty parameter
\(\alpha=10^4\).

Table~\ref{tab:mesh-refinement} reports the terminal-state and control errors
under simultaneous spatial and temporal mesh refinement. Both errors exhibit
approximately second-order convergence, with observed refinement orders
approaching \(2\). These results provide numerical evidence that the fully
discrete approximation accurately reproduces the exact penalized solution for
this benchmark.

\begin{table}[ht]
\centering
\caption{Mesh-refinement errors and observed convergence orders
for the fully discrete Crank--Nicolson approximation relative
to the exact penalized modal solution, with
\(y_T(x)=\sin(\pi x)\) and \(\alpha=10^4\).}
\label{tab:mesh-refinement}
\vspace{0.5em}
\begin{tabular}{cccccc}
\hline
$m$ & $n$
& Terminal error & Order
& Control error & Order \\
\hline
10 & 20
& $1.1397\times10^{-5}$ & --
& $5.2906\times10^{-2}$ & -- \\

20 & 40
& $2.8563\times10^{-6}$ & 1.996
& $1.3467\times10^{-2}$ & 1.974 \\

40 & 80
& $7.1452\times10^{-7}$ & 1.999
& $3.3821\times10^{-3}$ & 1.993 \\

80 & 160
& $1.7866\times10^{-7}$ & 2.000
& $8.4648\times10^{-4}$ & 1.998 \\
\hline
\end{tabular}
\end{table}

\FloatBarrier

\bibliographystyle{plain}
\bibliography{soft_to_hard_refs}

\end{document}